\documentclass[12pt,noamsfonts]{amsart}

\usepackage{amsfonts}
\usepackage{amssymb}
\usepackage{amscd}

\newcommand{\marginlabel}[1]%
  {\mbox{}\marginpar{\raggedleft\hspace{0pt}\bfseries\sf#1}}

\def\NN{{\mathbb N}}

\def\CC{{\mathbb C}}

\def\QQ{{\mathbb Q}}
\def\PP{{\mathbf P}}

\def\cI{\mathcal{I}}

\def\cA{\mathcal{A}}

\def\cO{\mathcal{O}}

\newcommand{\fra}{\frak{a}}

\DeclareMathOperator{\Sym}{Sym}

\DeclareMathOperator{\Cont}{Cont}

\DeclareMathOperator{\codim}{codim}

\DeclareMathOperator{\Spec}{Spec}

\DeclareMathOperator{\lc}{lc}

\DeclareMathOperator{\ord}{ord}

\newtheorem{lemma}{Lemma}[section]
\newtheorem{theorem}[lemma]{Theorem}
\newtheorem{corollary}[lemma]{Corollary}

\theoremstyle{definition}

\newtheorem{remark}[lemma]{Remark}

\newtheorem{example}[lemma]{Example}
\newtheorem{question}[lemma]{Question}

\theoremstyle{remark}
\newtheorem*{remark*}{Remark}
\newtheorem*{note*}{Note}

\begin{document}

\title{Multiplier ideals of hyperplane arrangements}

\author[M. Musta\c{t}\v{a}]{Mircea~Musta\c{t}\v{a}}
\address{Department of Mathematics, University of
Michigan\\ Ann Arbor, MI 48109, USA}  
\email{mmustata@umich.edu}

\thanks{The author
served as a Clay Mathematics Institute Research  Fellow
while this research has been done.}

\subjclass{Primary 14B05; Secondary 52C35}
\keywords{Arcs, multiplier ideals, hyperplane arrangements}

\begin{abstract}
In this note we compute multiplier ideals of hyperplane arrangements.
This is done using the interpretation of multiplier ideals in terms of
spaces of arcs from \cite{ELM}.
\end{abstract}

\maketitle

\bigskip

\section*{Introduction}

If $X$ is a complex smooth variety, and if $\fra\subseteq\cO_X$
is a coherent sheaf of ideals on $X$, then for every $\lambda\in\QQ_+^*$
one can define the multiplier ideal of exponent $\lambda$ denoted
$\cI(X,\fra^{\lambda})$. The interest in these ideals comes from the
fact that on one hand, they encode properties of the singularities
of the subscheme defined by $\fra$, and on the other hand, they appear
naturally in vanishing theorems.
This led to various applications in recent years, in higher dimensional
birational geometry as well as in local algebra (see \cite{lazarsfeld}
for a thorough presentation).

The definition of multiplier ideals is given in terms of a log
resolution of singularities for the pair $(X,\fra)$. This makes
the computation of specific examples quite delicate, so there are
not many classes of ideals for which the corresponding multiplier
ideals are known. One notable
exception is the case of multiplier ideals for
monomial ideals in the affine space, which have been 
computed by Howald using toric methods in \cite{howald}.
Most of the other examples where a general formula is known
proceed by reduction to the monomial case under suitable genericity
assumptions (see \cite{howald1}).

In this note we give a formula for the multiplier ideals of
hyperplane arrangements. This is achieved using the description
of multiplier ideals in terms of spaces of arcs from \cite{ELM}.

Let $\cA$ be a hyperplane arrangement in 
a vector space $V\simeq\CC^n$ with $n\geq 1$, i.e. $\cA$
is the (reduced) union of affine hyperplanes $H_1,\ldots,H_d\subset V$.
If $h_i$ is an equation defining $H_i$, then $\cA$
is defined by $f=\prod_ih_i$, and our goal is to compute $\cI(V, f^{\lambda})$.

We will assume that $\cA$ is central, i.e. that
all hyperplanes $H_i$ pass through the origin. This is not
an important restriction: if $\cA$ is an arbitrary arrangement,
and if we want to compute multiplier ideals in a neighbourhood of
a point $x\in V$, then it is enough to do the computation for
the subarrangement consisting of those hyperplanes passing through $x$.
Therefore after a suitable translation, we are in the case of a central arrangement.

Let $L(\cA)$ be the intersection lattice of $\cA$ (see \cite{OT}
for basic notions about arrangements). $L(\cA)$ is the
set of all intersections of the hyperplanes in $\cA$.
If $W\in L(\cA)$, then its rank is defined by $r(W):=\codim(W)$. 
We denote also by $s(W)$ the number of hyperplanes 
in $\cA$ containing $W$. Note that $V\in L(\cA)$
and we put $L'(\cA)=L(\cA)\setminus\{V\}$.
For every $W\in L(\cA)$, the
ideal defining $W$ in $V$ is denoted by $I_W$.
If $x$ is a real number, then $\lfloor x\rfloor$
denotes the largest integer $m$ such that $m\leq x$.
Our main result is the following

\begin{theorem}\label{main}
Let $\cA$ be a central hyperplane arrangement
defined by $(f=0)$. For every $\lambda\in\QQ_+^*$, we have
\begin{equation}\label{formula1}
\cI(V, f^{\lambda})=
\bigcap_{W\in L'(\cA)}I_W^{\lfloor \lambda s(W)\rfloor-r(W)+1}.
\end{equation}
\end{theorem}
Of course, we follow the convention that
if $r\leq 0$, then $(I_W)^r=\cO_V$.

Multiplier ideals are characterized by suitable conditions
on the orders of vanishing along divisors over $V$. One can interpret
the above theorem as saying that in this description it is enough to
consider only those divisors which appear as exceptional divisors for
the blowing-ups of $V$ along the various $W$ 
in $L'(\cA)$ (see Remark~\ref{valuations} below).

\begin{corollary}\label{cor0}
With the above notation, the support of the subscheme defined 
by $\cI(V,f^{\lambda})$ is equal to the union of those $W\in L'(\cA)$
such that $\lambda\geq\frac{r(W)}{s(W)}$.
\end{corollary}

An important invariant of an ideal $\fra$ is its log canonical threshold
$\lc(\fra)$.
It can be defined as the smallest $\lambda$ such that $\cI(X,\fra^{\lambda})
\neq\cO_X$.
We deduce the following formula in the case of a hyperplane arrangement.

\begin{corollary}\label{cor}
If $\cA$ is a central hyperplane arrangement defined by $(f=0)$, then
\begin{equation}
\lc(f)=\min_{W\in L'(\cA)}\frac{s(W)}{r(W)}.
\end{equation}
\end{corollary}

We will prove the above results in the next section, after we review 
the description of multiplier ideals in terms of spaces of arcs from 
\cite{ELM}. In the last section we discuss jumping coefficients 
in the sense of \cite{ELSV} for hyperplane arrangements.

\section{Proof of the main result}

We denote by $\NN$ the set of nonnegative integers.
Let us recall first a few basic facts about spaces of arcs.
For more details we refer to \cite{ELM}.
Recall that we work in the ambient smooth variety $V\simeq\CC^n$.
The arc space $V_{\infty}$ of $V$ is the set of all morphisms 
$\gamma : \Spec\CC[[t]]\longrightarrow V$. Since $V=\Spec(\Sym(V^*))$,
an arc $\gamma$ can be identified with a ring homomorphism
$\gamma^{\#} : \Sym(V^*)\longrightarrow \CC[[t]]$, hence with a morphism
of vector spaces $u : V^*\longrightarrow\CC[[t]]$. 
By looking at the coefficients of the different powers of $t$,
we may identify $u$ with $(u^{(m)})_m\in V^{\NN}$. Note 
also that if we fix a basis $e_1,\ldots,e_n$ for $V$ then
giving $u$ is the same as giving an element $(u_i)_i\in (\CC[[t]])^n$,
where $u_i=\gamma^{\#}(e_i^*)$. We will freely use these
identifications below.

If $g\in\cO(V)=\Sym(V^*)$ and if $\gamma\in V_{\infty}$
then the order $\ord_{\gamma}(g)$ 
of $g$ along $\gamma$ is defined as the order
of $\gamma^{\#}(g)\in\CC[[t]]$ (the order of zero is infinity).
If $m\in\NN$, then the contact locus $\Cont^{\geq m}(g)$ is the set
of all $\gamma$ such that $\ord_{\gamma}(g)\geq m$.

The space of $m$-jets of $V$ is defined as the set of all morphisms
$\delta : \Spec\,\CC[[t]]/(t^{m+1})\longrightarrow V$.
As above this can be naturally identified with $(\CC[[t]]/(t^{m+1}))^n$,
hence $V_m\simeq\CC^{(m+1)n}$. Every arc $\gamma\in V_{\infty}$ 
induces a jet $\gamma_m\in V_m$ by truncation, and this defines
a projection morphism $\psi_m : V_{\infty}\longrightarrow V_m$.
If a subset $C$ of $V_{\infty}$ is of the form $\psi_m^{-1}(S)$
for a Zariski closed subset $S$ of $V_m$ (such a set is called a 
\emph{closed cylinder}), 
then 
the codimension $codim(C)$ of $C$ in $V_{\infty}$ is the codimension
of $S$ in $V_m$ (it is easy to check that this is independent of $m$).
In particular, this applies to the above contact loci. 
If $C$ and $S$ are as above,
then $C$ is called irreducible if $S$ is so.
Note that every closed cylinder admits a unique decomposition
in irreducible components.

We recall now the description of multiplier ideals in terms of
contact loci. In fact, it is more convenient to change
the indexing of the multiplier ideals.
Recall that if $f\in\cO(V)$ and if $\lambda\in\QQ_+^*$, then
there is $\epsilon>0$ such that $\cI(V,f^{\lambda})=
\cI(V,f^t)$ for every $t\in [\lambda-\epsilon,\lambda)$. We
remark also that there are
only finitely many distinct multiplier ideals of exponent $\leq\lambda$
(for this and other basic properties of multiplier ideals, we refer to
\cite{lazarsfeld}).
We put
\begin{equation}
\widetilde{\cI}(V,f^{\lambda}):=\cI(V,f^{\lambda-\epsilon}),
\end{equation}
where $0<\epsilon\ll 1$. It follows that $\cI(V,f^{\lambda})
=\widetilde{\cI}(V, f^{\lambda+t})$ for $0<t\ll 1$.
With this notation, Corollary~3.15 in \cite{ELM} 
can be rewritten as

\begin{theorem}\label{description}
Let $f\in\cO(V)$ be a nonzero element and let $g\in\cO(V)$.
If $\lambda\in\QQ_+^*$, then $g\in\widetilde{\cI}(V,f^{\lambda})$
if and only if for every $m$ and every irreducible component $C$
of $\Cont^{\geq m}(f)$ we have
\begin{equation}
\ord_{\gamma}(g)\geq \lambda m-\codim(C)
\end{equation}
for all $\gamma\in C$.
\end{theorem}

We can give now the proof of our main result.

\begin{proof}[Proof of Theorem~\ref{main}]
Using the notation introduced above, proving formula~(\ref{formula1})
for all $\lambda$ is equivalent to proving that 
\begin{equation}\label{formula2}
\widetilde{\cI}(V,f^{\lambda})=\bigcap_{W\in L'(\cA)}I_W^{\lceil \lambda s(W)\rceil
-r(W)}
\end{equation}
for all $\lambda$, where $\lceil x\rceil$ denotes the smallest integer $m$
such that $x\leq m$.

Note first that in order to prove $(\ref{formula2})$, we may assume
that $\cA$ is an essential arrangement, i.e. that $\{0\}\in L(\cA)$. 
Indeed, let $T$ be the intersection of the hyperplanes in $\cA$. Then
there is an isomorphism $V\simeq T\times V'$ and an essential
hyperplane arrangement $\cA'$ in $V'$ defined by $f'$
 such that $\cA$ is the pull-back
of $\cA'$ via the canonical projection $\pi : T\times V'\longrightarrow V'$.
As 
$$\widetilde{\cI}(V,f^{\lambda})
=\pi^{-1}\widetilde{\cI}(V',(f')^{\lambda})$$ 
for all $\lambda$, it is clear that it is enough to prove 
(\ref{formula2}) for $\cA'$. We may therefore assume that $\cA$ is essential.

We use now Theorem~\ref{description} to describe
$\widetilde{\cI}(V,f^{\lambda})$. Recall that we denote by $h_i$
an equation of the hyperplane $H_i$, so $f=\prod_ih_i$.
It is clear that we have a decomposition
\begin{equation}\label{eq1}
{\rm Cont}^{\geq m}(f)=\bigcup_{\sum_ia_i=m}\bigcap_{i=1}^d{\rm Cont}^{\geq
a_i}(h_i).
\end{equation}

Consider $u=(u^{(j)})_{j\in\NN}\in V_{\infty}$, where 
$u^{(j)}\in V$ for all $j$.
Note that $u\in {\rm Cont}^{\geq a_i}(h_i)$ if and only if 
$u^{(j)}\in H_i$ for all $j$ such that $0\leq j\leq a_i-1$.
Hence $u\in\bigcap_i{\rm Cont}^{\geq a_i}(h_i)$ if and only if
for all $j$ we have
$$u^{(j)}\in\bigcap_{a_i>j}H_i.$$

For every $a=(a_1,\ldots,a_d)\in\NN^d$ and every $j\in\NN$, consider
$$W_j(a)=\bigcap_{a_i>j}H_i\in L(\cA).$$
We get in this way a sequence $W_{\bullet}(a)$ such that $W_j(a)
\subseteq W_{j+1}(a)$ for every $j$, and
$W_m(a)=V$ if $m\gg 0$. 

It follows from the above description that
$\bigcap_i{\rm Cont}^{\geq a_i}(h_i)$ is a closed irreducible cylinder of
codimension $\sum_{j\in\NN}\codim\,W_j(a)$. 
Theorem~\ref{description} and the decomposition
(\ref{eq1}) imply that given a polynomial $g\in\cO(V)$, we have
$g\in\widetilde{\cI}(V,f^{\lambda})$ if and only if for every
$a_1,\ldots,a_d\in\NN$ we have 
\begin{equation}\label{eq2}
\ord_u(g)\geq \lambda\cdot 
\sum_{i=1}^da_i-\sum_{j\in\NN}\codim\,W_j(a_{\bullet})
\end{equation}
for all $u=(u^{(j)})_{j\in\NN}\in V_{\infty}$ 
with $u^{(j)}\in W_j(a)$ for every $j$.

We reformulate now the above condition
starting this time with a sequence $W_{\bullet}$.
Suppose that $W_{\bullet}=(W_m)_{m\in\NN}$ is a sequence of
elements in $L(\cA)$ such that
$W_j\subseteq W_{j+1}$ for every $j$, and such that $W_m=V$
for $m\gg 0$.

Given $W_{\bullet}$
as above, let $Z(W_{\bullet})$ be the set of arcs $u=(u^{(m)})_{m\in\NN}\in
V_{\infty}$
such that $u^{(j)}\in W_j$ for every $j$. We see that 
$Z(W_{\bullet})$ is a closed irreducible cylinder, and
if
$u$ is a general element in $Z(W_{\bullet})$, then
\begin{equation}\label{formula_for_b}
\ord_u(h_i)=b_i(W_{\bullet}):=\min\{m\in\NN\vert W_m\not\subseteq H_i\}.
\end{equation}
Note that if $W_{\bullet}$ corresponds to $a\in\NN^d$, as before,
i.e. if $W_m=W_m(a)$ for all $m$, then 
$W_{a_i-1}\subseteq H_i$ by definition, so $b_i\geq a_i$ for all $i$. 

\smallskip

We claim that if we start with $W_{\bullet}$, and if $b=(b_1(W_{\bullet}),\ldots,
b_d(W_{\bullet}))$, where $b_i(W_{\bullet})$ is defined by (\ref{formula_for_b}),
then $W_j=W_j(b)$ for all $j$. 
Indeed, let us denote $W'_j=\bigcap_{b_i>j}H_i$. By definition, 
$j<b_i$ implies $W_j\subseteq H_i$, so $W_j\subseteq W'_j$.
For the reverse inequality, since $W_j\in L(\cA)$, it is enough
to show that if $W_j\subseteq H_k$ for some $k$, then $W'_j\subseteq H_k$.
But $W_j\subseteq H_k$ implies $b_k>j$ so $W'_j\subseteq H_k$ by definition. 
Therefore $W_j=W'_j$ for all $j$, which proves our claim.

We reinterpret now the condition in~(\ref{eq2}).
We claim that $g\in\widetilde{\cI}(V, f^{\lambda})$ if and only if for every
sequence $W_{\bullet}$ as above we have
\begin{equation}\label{eq3}
{\rm ord}_u(g)\geq\lambda\cdot\sum_{i=1}^db_i-\sum_{j\in\NN}\codim\,W_j
\end{equation}
for every $u\in Z(W_{\bullet})$. 

To see this, note first that by the above discussion,
if we start with $W_{\bullet}$ and define $b$ by (\ref{formula_for_b}),
then
$W_j=W_j(b)$ for all $j$, so
(\ref{eq3}) follows from (\ref{eq2}) for $b\in\NN^d$.
Conversely, if we start with $a\in\NN^d$, if $W_j=W_j(a)$
for every $j$ and if $b$ is given by (\ref{formula_for_b}), 
then (\ref{eq2}) follows from (\ref{eq3})
since $b_i\geq a_i$ for all $i$. 

\smallskip

We fix now a sequence $W_{\bullet}$ as above, and we rewrite
condition~(\ref{eq3}) for this sequence.
Let $r_j$ denote the codimension of $W_j$.
We may choose a basis $e_1,\ldots,e_n$ for $V$ such that for every $j$
we have $W_j=\langle e_{r_j+1},\ldots,e_n\rangle$.

In this system of coordinates, we identify an arc $\gamma\in V_{\infty}$
with an element $u=(u_1,\ldots,u_n)\in(\CC[[t]])^n$, where $u_i=\gamma^{\#}
(X_i)$. Here $X_i=e_i^*\in\cO(V)=\Sym(V^*)=\CC[X_1,\ldots,X_n]$.
Therefore
$Z(W_{\bullet})$ can be described as
the set of those arcs $u=(u_i)_i$ such that for all $i$ and $j$
with $r_j\geq i$, we have $\ord_u(X_i)>j$. Hence $u\in V_{\infty}$
is in $Z(W_{\bullet})$ if and only if $\ord_u(X_i)\geq \min\{j\vert
r_j<i\}$. It follows that in this system of coordinates, the set
of those
$g\in\cO(V)$ which satisfy
(\ref{eq3}) for $W_{\bullet}$ 
is the monomial ideal generated by those $X_1^{\alpha_1}\ldots X_n^{\alpha_n}$
such that 
\begin{equation}\label{eq4}
\sum_{k=1}^n\alpha_k\cdot\min\{j\vert r_j<k\}
\geq\lambda\cdot\sum_{i=1}^db_i-\sum_{j\in\NN}r_j.
\end{equation}

Let $p_i:=\min\{j\vert r_j<i\}$ for every $i$ with $1\leq i\leq n$,
so we have $p_1\geq p_2\geq\ldots\geq p_n$.
We see that $W_{\bullet}$ consists of $p_n$ zero subspaces,
and of $p_{n-i}-p_{n-i+1}$ subspaces of dimension $i$, for every $1\leq i\leq n-1$.
Hence $\sum_{j\in\NN}r_j=\sum_{k=1}^np_k$,
and (\ref{eq4}) becomes
\begin{equation}\label{eq5}
\sum_{k=1}^n(\alpha_k+1)p_k\geq\lambda\cdot\sum_{i=1}^d b_i.
\end{equation}

We put $q_n=p_n$
and $q_i=p_i-p_{i+1}$ for $1\leq i\leq n-1$.
We call
the sequence $W_{\bullet}$ \emph{strict} if $\dim(W_j)=j$ for $0\leq j\leq n$.
It is clear that every $W_{\bullet}$ as above is obtained from a strict
$W'_{\bullet}$ by taking $q_{n-i}$ copies of each $W'_i$, 
for $0\leq i\leq n-1$
(the remaining $W_i$ 
being equal to the ambient space). Note that some of the 
$q_i$ may be zero.
In addition, if $b'=(b'_1,\ldots,b'_d)$ corresponds to $W'_{\bullet}$, 
then we have $b_i=p_{n+1-b'_i}$ for every $i$.
If for every $k$ with $1\leq k\leq n$ 
we denote by $\tau_k$ the number of elements in
 $\{i\vert b'_i=k\}$, we get
$$\sum_{i=1}^db_i=\sum_{k=1}^n\tau_{n+1-k}p_k.$$

We let now $W_{\bullet}$ vary for a fixed $W'_{\bullet}$,
i.e. we allow the $q_i$ to vary. It follows from (\ref{eq5})
that the set of those $g\in\cO(V)$ which satisfy (\ref{eq3})
for all such $W_{\bullet}$ is the monomial ideal generated by
those $X_1^{\alpha_1}\ldots X_n^{\alpha_n}$ such that
\begin{equation}\label{eq6}
\sum_{k=1}^n(\alpha_k+1)\cdot \sum_{j=k}^nq_j\geq\lambda\cdot\sum_{k=1}^n
\tau_{n+1-k}\cdot\sum_{j=k}^nq_j
\end{equation}
for every $q_1,\ldots,q_n\in\NN$.
Equivalently, we have
\begin{equation}\label{eq7}
\sum_{\ell=1}^k(\alpha_{\ell}+1)
\geq\lambda\cdot\sum_{\ell=1}^k\tau_{n+1-\ell}
\end{equation}
for every $k$, $1\leq k\leq n$.

We can reformulate this condition without any reference to
the system of coordinates as follows. If $g\in\cO(V)$, then
$g$ satisfies (\ref{eq3}) for all $W_{\bullet}$ corresponding
to a fixed strict sequence $W'_{\bullet}$ if and only if
for every $k$ with $1\leq k\leq n$ we have
\begin{equation}\label{eq8}
g\in I_{W_{n-k}}^{\lceil\lambda\cdot s_k\rceil-k}
\end{equation}
for every $k$, where $s_k=\sum_{\ell=1}^k\tau_{n+1-\ell}=s(W'_{n-k})$.
Since every $W\in L(\cA)$, $W\neq V$ appears as $W'_{n-k}$
for a suitable strict sequence $W'_{\bullet}$ (where $k=r(W)$),
we deduce that $g\in\widetilde{\cI}(V,f^{\lambda})$
if and only if for every $W\in L'(\cA)$ we have
$$g\in I_W^{\lceil \lambda s(W)\rceil -r(W)}.$$
This completes the proof of the theorem.
\end{proof}

\begin{remark}\label{valuations}
The multiplier ideals of $f$ 
are characterized by the orders of vanishing along divisors over $V$.
More precisely, suppose that $E$ is a divisor on a smooth variety
$V'$ such that there is a proper, birational morphism $\pi\colon
V'\to V$. Associated to such $E$ there are two numbers: the order
$a_E(f)$ of $E$ in the pull-back of ${\rm div}(f)$ to $V'$, and
the order $k_E$ of $E$ in the relative canonical divisor $K_{V'/V}$.
The multiplier ideal
${\mathcal I}(f^{\lambda})$ is the set of all $g\in {\mathcal O}(V)$
such that
$$a_E(g)>\lambda a_E(f)-k_E-1$$
for all $E$ as above.

We may reinterpret the statement of Theorem~\ref{main}
as saying that if $f$ defines a hyperplane arrangement ${\mathcal A}$,
then in the above description of ${\mathcal I}(f^{\lambda})$
it is enough to consider only those $E$ which are the exceptional
divisors of $V$ along the various $W$ in $L'({\mathcal A})$.
Indeed, for such $E$ corresponding to $W$, note that 
$a_E(f)=s(W)$ and $k_E=r(W)-1$. Moreover, we have $g\in I_W^m$
if and only if $a_E(g)\geq m$.
\end{remark}

\begin{proof}[Proof of Corollaries~\ref{cor0} and~\ref{cor}]
The description of the support of the subscheme defined by
$\cI(V,f^{\lambda})$ follows immediately from the formula
in Theorem~\ref{main}. We deduce that $\cI(V,f^{\lambda})=\cO_X$
if and only if $\lambda<r(W)/s(W)$ for every $W\in L(\cA)$. 
This gives the formula for $\lc(f)$.
\end{proof}

\section{Jumping coefficients of hyperplane arrangements}

We have seen that Theorem~\ref{main} gives a formula for the 
log canonical threshold of a hyperplane arrangement. In this section
we discuss the higher jumping coefficients
of multiplier ideals.

Let us recall the definition from \cite{ELSV}. If $X$ is a smooth
variety and if $\fra\subset\cO_X$ is a nonzero coherent sheaf of ideals, then
the jumping coefficients of $\fra$ are the elements of
$$\{\lambda\in\QQ_+^*\vert \cI(X,\fra^{\lambda})\neq\cI(X,\fra^{\lambda-
\epsilon})\,{\rm if}\,0<\epsilon\ll 1\}.$$
Hence $\lc(\fra)$ is the smallest jumping coefficient.

Suppose now that $X=V\simeq\CC^n$ and that $\fra$ is a principal ideal generated
by the nonzero polynomial $f$. We recall two basic facts
about the jumping coefficients of $f$ (see \cite{ELSV} for proofs).
First, $1$ is always a jumping coefficient.
Second, if $\lambda>1$, then $\lambda$ is a jumping coefficient if
and only if so is $\lambda-1$. Therefore in order to understand these
invariants it is enough to study those contained in $(0,1)$. 

The jumping coefficients are related also to the Bernstein polynomial 
$b_f\in\CC[s]$. 
More precisely, if $\lambda$ is a jumping coefficient in $(0,1]$,
then $b_f(-\lambda)=0$.

The following is an immediate consequence of Theorem~\ref{main}.

\begin{corollary}\label{cor2}
If $\cA$ is a central hyperplane arrangement defined by $(f=0)$,
then $\lambda\in\QQ_+^*$ is a jumping coefficient of $f$ if and only if
there are
$W\in L(\cA)$ and $m\in\NN$ such that $\lambda=\frac{r(W)+m}{s(W)}$
and such that
\begin{equation}
\bigcap_{W'\in L'(\cA)}I_{W'}^{\lceil \lambda s(W')\rfloor -r(W')}
\not\subseteq I_W^{m+1}.
\end{equation}
In fact, we may take the above intersection over only those
$W'\in L(\cA)$ containing $W$. 
\end{corollary}

Note that by the above Corollary, each $W\in L'(\cA)$
gives some candidates for jumping coefficients (we recall that
we consider
only those numbers in $(0,1)$).
However, unlike in the case of the log canonical threshold,
it is not at all clear how to interpret the
condition that $\lambda$ is a jumping number just in terms of the
intersection lattice $L(\cA)$.

\begin{question}
Given $\lambda\in\QQ_+^*$, does the condition that $\lambda$
is a jumping number for $f$ depend only on the intersection lattice 
$L(\cA)$ ? 
\end{question}

We consider now several particular cases.

\begin{example}\label{ex1}
We say that $\lambda$ is a set-theoretic jumping coefficient
for $f$ if the support of the subscheme defined by $\cI(V,f^{\lambda})$
strictly contains the support of the subscheme defined by
$\cI(V,f^{\lambda-\epsilon})$ for $0<\epsilon\ll 1$.
The smallest such invariant is the log canonical threshold.

As in the case of the log canonical threshold, these numbers are easy to describe:
$\lambda$ is a set-theoretic jumping coefficient of $f$
if and only if there is $W\in L'(\cA)$ such that $\lambda=\frac{r(W)}{s(W)}$
and such that $\frac{r(W)}{s(W)}\leq\frac{r(W')}{s(W')}$ for every
$W'\in L'(\cA)$ containing $W$. This follows immediately from the
description in Corollary~\ref{cor0}
of the support of the subscheme defined by 
$\cI(V,f^{\lambda})$. 
\end{example}

\begin{example}\label{ex2}
Suppose that $\cA$ is a generic arrangement of $d\geq n$ hyperplanes
in $\CC^n$. 
This means that for every $k\leq n$, the intersection of
any $k$ of the hyperplanes in $\cA$ has dimension $n-k$.
Let $f$ be an equation for $\cA$.
Walther proved in \cite{walther} that the roots of the Bernstein polynomial
$b_f(s)$ of $f$ are given by $-\frac{n+i}{d}$, with $0\leq i\leq 2d-2-n$.
We see that the negatives of the roots in $[-1,0)$ are precisely the
jumping numbers in $(0,1]$.
Indeed, if ${\mathbf m}=(X_1,\ldots,X_n)$ and if $\lambda\leq 1$, then it 
follows from Theorem~\ref{main} that
$$\cI(\CC^n,f^{\lambda})={\mathbf m}^{\lfloor \lambda d\rfloor-n+1}.$$
\end{example}

\begin{example}\label{ex3}
Let $\cA$ be a central arrangement,
and consider $W\in L(\cA)$ with $r(W)=2$. For every $W'\in\ L'(\cA)$
strictly containing $W$, $W'$ has to be one of the hyperplanes,
so $r(W')=s(W')=1$. It follows from Corollary~\ref{cor2} that
all candidates for jumping coefficients $\leq 1$ corresponding to $W$
$$\frac{2+m}{s(W)}\,\, {\rm for}\,\,0\leq m\leq s(W)-2$$
are indeed jumping coefficients. 
Of course, this is interesting only if $s(W)\geq 3$.
\end{example}

\begin{example}\label{ex4}
Suppose now that $W\in L(\cA)$ has codimension
three, and let us analyse
 when some of the candidates for jumping coefficients corresponding to $W$
are indeed jumping coefficients. It follows from Corollary~\ref{cor2}
that in order to decide this it is enough to consider in 
$V/W$ the arrangement induced by those hyperplanes in $\cA$ which
contain $W$. We may therefore assume that $\dim\,V=3$ and that $W=(0)$.
Since $\cA$ is a central arrangement, it is convenient to consider also
the induced arrangement $\overline{\cA}$ in $\PP(V)\simeq\PP^2$. 
If $P\in\overline{\cA}$ is a point in $L(\overline{\cA})$ we still
denote by $s(P)$ the number of lines in $\overline{\cA}$ passing through $P$.

Let ${\mathcal T}_0(\cA)$ be the set of jumping coefficients of 
$\cA$ in $(0,1)$ corresponding to codimension two elements in $L(\cA)$
\begin{equation}
{\mathcal T}_0(\cA):=\{\frac{j}{s(W)}\vert 2\leq j\leq s(W)-1,\,r(W)=2\}. 
\end{equation}

Let $\lambda\in (0,1)$ be a candidate corresponding to the zero subspace,
hence $\lambda=j/d$, where $d$ is the number of hyperplanes in $\cA$
 and $3\leq j\leq d-1$. 
We may assume that $\lambda\not\in {\mathcal T}_0(\cA)$, and we want
to determine whether $\lambda$ is a jumping coefficient.
It follows from Corollary~\ref{cor2} that
$\lambda$ is a jumping coefficient if and only if
there is a hypersurface in $\PP(V)$ of degree $\leq (j-3)$
 passing through each point $P$ in $L(\overline{\cA})$
with multiplicity $\geq \lceil js(P)/d\rceil-2$.

Suppose first that $j=3$. We see that $3/d$ is a jumping coefficient
if and only if for every point $P\in L(\overline{\cA})$ 
we have $s(P)\leq 2d/3$.

Suppose now that $j=4$. It follows from the above discussion that
$4/d$ is a jumping coefficient if and only if there is no point
$P\in L(\overline{\cA})$ with $s(P)>3d/4$, and if there are no three 
noncolinear
points $Q_1,Q_2,Q_3$ with $s(Q_i)>d/2$ for all $i$.
However, it is easy to see that if $s(Q_1)$, $s(Q_2)>d/2$, then
$\overline{\cA}$ 
consists precisely of $(d+1)/2$ lines passing through $Q_1$ and
$(d+1)/2$ lines passing through $Q_2$, and the line joining $Q_1$ and $Q_2$
belongs to $\overline{\cA}$. In particular, $s(P)=2$ for every other point 
$P\in L(\overline{\cA})$. We conclude that $4/d$ is a jumping coefficient
of $\cA$ if and only if for every point $P\in L(\cA)$ we have 
$s(P)\leq 3d/4$.

Consider now the case $j=5$.
We see that $5/d$ is a jumping coefficient if and only
if there is a conic passing through each point $P\in L(\overline{\cA})$
with multiplicity $\geq \lceil 5s(P)/d\rceil -2$. In order for this to be true,
we clearly need that $s(P)\leq 4d/5$ for every point $P\in L(\overline{\cA})$.
We show that this is also sufficient. Note first that we can't have
two points $P_1$, $P_2\in L(\overline{\cA})$ with $s(P_1)$, $s(P_2)>3d/5$
or four points $Q_1,\ldots,Q_4$ with $s(Q_i)>2d/5$ for all $i$. 
Therefore the only case we need to consider is when we have two points
$P_1$, $P_2\in L(\overline{\cA})$ with $s(P_1)>3d/5$ and $s(P_2)>2d/5$.
Again, it is easy to see that in this case every line passes through
either $P_1$ or $P_2$, and the line $\ell$ joining $P_1$ and $P_2$
is in $\overline{\cA}$. In particular,
$s(Q)=2$ for every other point $Q$.
In this case we may take the conic to be the union
of $\ell$ and another line in 
$\cA$ passing through $P_1$. 
We conclude that $5/d$ is a jumping coefficient if and only if 
$s(P)\leq 4d/5$ for every point $P\in L(\overline{\cA})$.

We consider now the case 
$j=d-1$. We see that $(d-1)/d$ is a jumping coefficient
if and only if there is a hypersurface of degree $\leq (d-4)$ 
passing through each point $P\in L(\overline{\cA})$ with multiplicity
$\geq s(P)-2$ (we use the fact that $s(P)<d$ for every such $P$).
It is clear that this is not possible if there is $P\in L(\overline{\cA})$
such that $s(P)=d-1$. 
Conversely, if this is not the case, then it is easy to see
that we can find four lines in $\overline{\cA}$ no three of them meeting
in one point. We may take the hypersurface as the union of the remaining
$d-4$ lines. Hence $(d-1)/d$ is a jumping coefficient
if and only if there is no $P\in L(\overline{\cA})$
with $s(P)=d-1$.
\end{example}

\begin{question}
With the notation in the above example, let $\lambda=\frac{j}{d}$
with $3\leq j\leq d-1$, and suppose that $\lambda\not\in {\mathcal T}_0(\cA)$.
Is it true that $\lambda$ is a jumping number if and only if for every
point $P\in\overline{\cA}$ we have $s(P)\leq\frac{(j-1)d}{j}$ ?
The above example shows that this is the case if $j\in\{3,4,5,d-1\}$.
Note also that in each of these cases, we can choose the relevant
hypersurface as a union of lines in $\overline{\cA}$.
\end{question}

\subsection*{Acknowledgements.}
We are grateful to Zach Teitler for his comments
on a previous version of this paper, and to the referee 
for pointing out Remark~\ref{valuations}.

\providecommand{\bysame}{\leavevmode \hbox \o3em
{\hrulefill}\thinspace}


\begin{thebibliography}{ELM}


\bibitem[ELM]{ELM}
L.~Ein, R.~Lazarsfeld and M.~Musta\c{t}\v{a},
Contact loci in arc spaces, Compos. Math. \textbf{140} (2004), 1229--1244.

\bibitem[ELSV]{ELSV}
L.~Ein, R.~Lazarsfeld, K.~E.~Smith and D.~Varolin,
Jumping coefficients of multiplier ideals,
Duke Math. J. \textbf{123} (2004), 469--506.


\bibitem[How1]{howald1}
J.~Howald, Multiplier ideals of sufficiently general polynomials,
preprint, math.AG/0303203.

\bibitem[How2]{howald}
J.~Howald, Multiplier ideals of monomial ideals, Trans. Amer. Math. Soc.
\textbf{353} (2001), 2665--2671.


\bibitem[Laz]{lazarsfeld}
R.~Lazarsfeld, \emph{Positivity in algebraic geometry II}, 
Ergebnisse der Mathematik und ihrer Grenzgebiete. 3. Folge, 
A series of Modern Surveys in Mathematics, Vol. 49,
Springer-Verlag, Berlin, 2004.



\bibitem[OT]{OT}
P.~Orlik and H.~Terao, \emph{Arrangements of hyperplanes},
Grundlehren der mathematischen Wissenschaften, Vol. 300, Springer-Verlag,
New York, 1992.

\bibitem[Wal]{walther}
U.~Walther, Bernstein-Sato polynomial versus cohomology
of the Milnor fiber for generic hyperplane arrangements,
Compos. Math., to appear.

\end{thebibliography}
\end{document}